\documentstyle[12pt]{article}
\def\ra{{\rightarrow}}

\def\R{{R}}

\title {\bf Stochastic differential equations with non-Lipschitz coefficients:
I. Pathwise uniqueness and Large deviations}
\begin{document}
\author{Shizan FANG \ and Tusheng ZHANG}
\date{}
\maketitle
\medskip
 I.M.B, UFR Sciences et techniques, Universit\'e de Bourgogne, 9 avenue Alain Savary, B.P. 47870, 21078 Dijon, France.\\
Department of Mathematics, University of Manchester, Oxford road, Manchester, M13 9PL, England.
\vskip 10mm
\abstract{We study a class of  stochastic differential equations with non-Lipschitzian coefficients. A unique strong solution is obtained and a large deviation principle of Freidlin-Wentzell type has been  established.}
\bigskip
\vskip 15mm
\section{Introduction}
\bigskip

\noindent Let $ \sigma:\ R^d\rightarrow  R^d\otimes R^m$ and 
$b:\ R^d\rightarrow R^d$ be continuous functions. It is well-known that the following It\^o s.d.e:

\begin{equation}
 dX(t)=\sigma(X(t))\,dW_t + b(X(t))\, dt,\quad X(0)=x_o
\end{equation}
has a weak solution up to a lifetime $\zeta$ (see [SV], [IW, p.155-163]), where $t\rightarrow  W_t$ is a $\R^m$-valued standard Brownian motion. It is also known that under the assumption of linear growth of coefficients $\sigma$ and $b$, the lifetime $\zeta= +\infty$ almost surely. If the s.d.e $(1)$ has the pathwise uniqueness, then it admits a strong solution (see [IW, p.149], [RY, p.341]). So the study of pathwise uniqueness is of great interest. It is a classical result that under the Lipschitz conditions, the pathwise uniqueness holds and the solution of s.d.e. $(1)$ can be constructed using Picard interation; morever the solution depends on the  initial values continuousely. The main tool to these studies is the Gronwall lemma. When the coefficients $\sigma$ and $b$ do not satisfy the Lipschitz conditions, the use of Gronwall lemma meets a serious difficulty. Therefore, there are very few results of pathwise uniqueness of  solutions of s.d.e. beyond the Lipschitzian (or locally Lipschitzian ) conditions in the literature except in the one dimensional case (see [IW, p.168], [RY, CH IX-3]). In the case of ordinary differential equations, the Gronwall lemma was generalized in order to establish the uniqueness result (see e.g. [La]). However the method is not applicable to s.d.e.. In this work, we shall deal with a class of non-Lipschitzian s.d.e.. Namely, we shall assume that

$$ \left\{\matrix{||\sigma(x)-\sigma(y)||^2&\leq& C\, |x-y|^2\,\log{{ 1\over
 |x-y|}},
&\quad\hbox{\rm for  } |x-y|<1,\cr
|b(x)-b(y)|&\leq& C\, |x-y|\,\log{{ 1\over |x-y|}},&\quad\hbox{\rm for  }
 |x-y|<1}\right.
\leqno(H1)$$
where $|\cdot|$ denotes the Euclidean distance in $\R^d$ and 
$ ||\sigma||^2=\sum_{i=1}^d\sum_{j=1}^m \sigma_{ij}^2$. 
\vskip 0.3cm

We will prove the pathwise uniqueness of solutions under $(H.1)$ and non explosion of the solution under the growth condition $|x|\log |x|$. The  results are valid for any dimension.  Our idea  is to construct a family of positive increasing functions $(\Phi_{\rho})_{\rho >0}$ on $R_+$ so that the Gronwall lemma can be applied to the composition of the functions $(\Phi_{\rho})$ with appropriate processes. This family of positive functions   plays a crucial role. In this work, we will also establish a Freidlin-Wentzell type large deviation principle for the solutions of the s.d.e's (see [FW]). As a by-product, it is seen that the  unique solution of the s.d.e.  can be obtained by Euler approximation. Our strategy for the large deviation is to prove that the Euler approximations to the s.d.e. is exponentially fast. The method of estimating moments used in the literature ([DS],[DZ], [S] ) wouldn't work here because of the non-Lipschitzian feature of the coefficients. We again appeal to a family of positive functions  $(\Phi_{\rho, \lambda })_{\rho >0}$. The proof of the uniform convergence of solutions of the corresponding skeleton equations over compact level sets is also tricky due to the non-Lipschitzian feature. 
\medskip

\quad The organization of the paper is as follows. In section 2, we shall discuss the case of ordinary differential equations; although the results about non explosion and uniqueness are not new, but our method can also be used to study the dependence of initial values and the non confluence of the equations.  In section 3, we shall consider the s.d.e. The pathwise uniqueness and the criterion of non-explosion will be established. However, the supplementary difficulties will appear when we deal with the dependence with respect to initial values and the non confluence of s.d.e., that we shall study in a forthcoming paper. The ordinary differential equation

\begin{equation}
 {dX(t)\over dt}=b(X(t)),\quad X(0)=x_o
\end{equation}
gives rise to a dynamical system on $\R^d$. In section 4, we shall consider its small perturbation by a white noise. Namely, we shall consider the s.d.e

\begin{equation}
 dX^{\varepsilon}(t) = \varepsilon^{1\over 2} \, \sigma(X^{\varepsilon}(t))\, dW_t + b(X^{\varepsilon}(t) )\, dt,\quad X^{\varepsilon}(0) = x_o
\end{equation}
and state  a large deviation principle for $ (X^{\varepsilon}(t))_{t\in [0,1]}$. Section 5 and 6 are devoted to the proof of the large deviation principle. Section 5 is for the case of bounded coefficients. Section 6 is for the general case.

\vskip 10mm
\section{Ordinary differential equations}
\bigskip

\noindent  Let $b: \R^d\rightarrow  \R^d$ be a continuous function. It is essentially due to Ascoli-Arzela theorem that the differential equation $(2)$ has a solution up to a lifetime $\zeta$. The following result weakens  the linear growth condition for non explosion.

\medskip
\noindent {\bf Theorem 2.1}\quad {\it Let $r: \R_+\rightarrow \R_+$ be a continuous function such that
\medskip
$(i)$\quad $\lim_{s\rightarrow +\infty}r(s)=+\infty$,
\smallskip
$(ii)$\qquad $\int_0^\infty {ds\over sr(s)+1}=+\infty$.
\smallskip

Assume that it holds
\begin{equation}
 |b(x)|\leq C\ \bigl(|x|r({|x|^2}) + 1\bigr).
\end{equation}
Then the lifetime is infinte: $\zeta = +\infty$.}

\medskip
\noindent{\bf Proof.}\quad Define for $\xi\geq 0$,
$$ \psi(\xi)=\int_0^\xi {ds\over sr(s)+1}
\quad\hbox{\rm and  } \Phi(\xi) = e^{\psi(\xi)}.$$
We have
\begin{equation}
\Phi'(\xi) = {\Phi(\xi)\over \xi r(\xi)+1}.
\end{equation}
Let  $\xi_t=|X_t|^2$, where $X(t)$ is a solution to (2).  Then

\begin{equation}
 {d\over dt}\Phi(\xi_t) = 2\Phi'(\xi_t)\bigl<X_t, b(X(t))\bigr>,
\end{equation}
where $\bigl<\ ,\ \bigr>$ denotes the inner product in $\R^d$. By assumption $(4)$, we have
\begin{equation}
\Bigl| {d\over dt}\Phi(\xi_t)\Bigr|
\leq 2C \Phi'(\xi_t)|X_t|\,\bigl(|X_t|r(\xi_t)+1\bigr).
\end{equation}
By $(i)$, it holds that
$$  \sup_{s\geq 0}{s^2 r(s^2) + s\over s^2r(s^2)+ 1 } <+\infty.$$

\noindent Therefore for some constant $C_2>0$,
\begin{equation}
\Bigl| {d\over dt}\Phi(\xi_t)\Bigr|
\leq 2C_2 \Phi(\xi_t).
\end{equation}

It follows that for $t<\zeta$,
$$ \Phi(\xi_t)\leq \Phi(|x_o|^2) + 2C_2\int_0^t \Phi(\xi_s)\, ds.$$
By Gronwall lemma, we have
\begin{equation}
 \Phi(\xi_t) \leq \Phi(|x_o|^2)\,e^{2C_2\,t}.
\end{equation}
If $\zeta <+\infty$, letting $t\uparrow \zeta$ in $(9)$, we get 
$\Phi(\xi_\zeta)\leq \Phi(|x_o|^2)\,e^{2C_2\zeta}$ which is impossible because of $\xi_\zeta=+\infty$, $\Phi (+\infty)=+\infty$.

\medskip
\noindent{\bf Remark.}\quad By considering the inequality
${d\over dt}\Phi(\xi_t)\geq -2C_2\, \Phi(\xi_t)$, we have
$$\Phi(\xi_t)\geq \Phi(|x_o|)-2C_2\int_0^t \Phi(\xi_s)\, ds,$$
which yields to
$$\Phi(\xi_t)\geq \Phi(|x_o|)e^{-2C_2 t}.$$
If we denote by $X_t(x_o)$ the solution to $(2)$ with initial value $x_o$, then we get
$\lim_{|x_o|\ra +\infty}\Phi(|X_t(x_o)|)=+\infty$,
which implies that 
\begin{equation}
\lim_{|x_o|\ra +\infty}|X_t(x_o)|=+\infty.
\end{equation}

\noindent  In what follows, to be simplified, we shall assume that the solutions of $(2)$ 
have non explosion.
\medskip

\noindent {\bf Theorem 2.2}\quad {\it Let $r: ]0,1[\rightarrow  \R_+$ be a continuous function such that
\smallskip

\quad (i)\qquad $ \lim_{s\rightarrow  0} r(s)=+\infty$;

\quad (ii)\qquad $ \int_0^a {ds\over sr(s)}=+\infty$ for any $a>0$.

Assume that
\begin{equation}
  |b(x)-b(y)|\leq C\, |x-y|\ r(|x-y|^2)\quad\hbox{\rm for  } |x-y|<1.
\end{equation} 

Then the differential equation $(2)$ has an unique solution.}
\medskip

\noindent {\bf Proof.}\quad Let $ (X(t))_{t\geq 0}$ and $ (Y(t))_{t\geq 0}$ be two solutions of the equation $(2)$. Set $\eta_t = X(t)-Y(t)$ and $ \xi_t=|\eta_t|^2$. Let $\rho>0$, consider
$$ \psi_\rho(\xi)=\int_0^\xi {ds\over sr(s)+\rho}
\quad\hbox{\rm and}\quad \Phi_\rho(\xi)=e^{\psi_\rho(\xi)}.$$
We have
\begin{equation}
\Phi_\rho'(\xi) = {\Phi_\rho(\xi)\over \xi\ r(\xi)+\rho}.
\end{equation}
Let
$$ \tau=\inf\{ t>0,\ \xi_t\geq 1/2\,\}.$$
By assumption $(11)$, we have
$$\bigl|\bigl<\eta_t, b(X(t))-b(Y(t))\bigr>\bigr|
\leq C\ \xi_t\,r(\xi_t),\quad t<\tau.$$

Therefore according to $(12)$, for $t<\tau$, by chain rule,
$$ \Phi_\rho(\xi_t) \leq 1 + 2C\, \int_0^t \Phi_\rho(\xi_s)\, ds,$$
which implies that
$\Phi_\rho(\xi_t)\leq e^{2C\,t}$ for $t<\tau$.
Letting $\rho\downarrow 0$, we get that $ e^{\psi_0(\xi_t)}\leq e^{2Ct}$. Now by hypothesis 
(ii), we obtain that $\xi_t=0$ for all $t<\tau$. If $\tau<+\infty$, letting $t\uparrow\tau$, we get
$$ {1\over 2} = \xi_\tau = 0,$$
which is absurd. Therefore $\xi_t=0$ for all $t\geq 0$. In other words, 
$X(t)=Y(t)$ for $t\geq 0$.

\medskip
\noindent {\bf Example 2.3}\quad Define
$$ f(x_1,x_2) = \sum_{k\geq 1} {\sin{kx_1}\cdot\sin{k x_2}\over k^2}.$$
Obviously the function $f$ is continuous on $\R^2$. We have
\begin{equation}
 |f(X)-f(Y)|\leq C\ |X-Y|\log{{1\over |X-Y|}}\quad\hbox{\rm for  } |X-Y|< {1\over e}
\end{equation}
where $X=(x_1,x_2)$ and $Y=(y_1,y_2)$. In fact,

$$ f(X)-f(Y)=\sum_{k=1}^\infty \Bigl\{{(\sin{kx_1}-\sin{ky_1})\sin{kx_2}\over k^2}
+  {(\sin{kx_2}-\sin{ky_2})\sin{ky_1}\over k^2}\Bigr\}.$$
It follows that

$$ |f(X)-f(Y)|\leq 2\sum_{k=1}^\infty \Bigl\{{|\sin{(k{x_1-y_1\over 2})}|\over k^2}
+ {|\sin{(k{x_2-y_2\over 2})}|\over k^2}\Bigr\}.$$

\noindent {\bf Lemma 2.4}\quad {\it For $ 0<\theta< 1/e$, we have
\begin{equation}
 V(\theta):= \sum_{k=1}^\infty {|\sin{k\theta}|\over k^2}
\leq C_1\ \theta\,\log{{1\over\theta}}.
\end{equation}}
\medskip

\noindent {\bf Proof.}\quad Consider $ \phi(s)={\sin{s\theta}\over s^2}$. We have
$$ \phi'(s) = {s^2\theta\cos{s\theta}-2s\sin{s\theta}\over s^4}.$$
Then $|\phi'(s)|\leq {3\theta\over s^2}$. Let 
$ W(\theta)=\int_1^{+\infty} {|\sin{s\theta}|\over s^2}\ ds$. We have
\begin{eqnarray}
 |V(\theta)-W(\theta)|
&\leq & \sum_{k=1}^{+\infty} \int_k^{k+1} |\phi(s)-\phi(k)|\ ds\nonumber \\
&\leq & 3\theta\sum_{k=1}^{+\infty} {1\over k^2}= {\pi^2\theta\over 2}.
\end{eqnarray}
Now
$$ W(\theta) = \theta\,\int_\theta^{+\infty} {|\sin{t}|\over t^2}\, dt
\leq \theta\, \int_{\theta}^1 {\sin{t}\over t}{dt\over t}
+ \theta\,\int_1^{+\infty}{ds\over s^2}$$
which is dominated by
$$ \theta\Bigl(\log{{1\over\theta}}+1\Bigr).$$
Therefore, according to $(15)$
$$ V(\theta) \leq\theta\Bigl(\log{{1\over\theta}}+1+{\pi^2\over 2}\Bigr),$$
which is less that
$2({\pi^2\over 2}+1)\ \theta\log{{1\over\theta}}$ for $0<\theta<{1\over e}$.

\medskip
\noindent Now applying $(14)$, for $|x_1-y_1| + |x_2-y_2| < 1/e$,
\begin{eqnarray}
 &|f(X)-f(Y)|\nonumber\\
\leq & 4C_1 \Bigl( {|x_1-y_1|\over 2}\log{{1\over |x_1\!-\!y_1|}}
+{|x_2-y_2|\over 2}\log{{1\over |x_2\!-\!y_2|}}\Bigr)\nonumber \\
\leq & 4C_1\Bigl({|x_1-y_1|+|x_2-y_2|\over 2}\log{{2\over |x_1\!-\!y_1|+|x_2\!-\!y_2|}} \Bigr)
\end{eqnarray}
where the last inequality was due to the concavity of the function $\xi\log{{1\over\xi}}$ over
$]0,1[$. Therefore $(13)$ holds for some constant $C>0$.

\medskip
\noindent In what follows, we shall study the dependence of initial values.
\medskip

\noindent {\bf Theorem 2.5}\quad {\it Assume that the conditions $(4)$ and $(11)$ hold. Then
$ x_o\ra X_t(x_o)$ is continuous.}

\medskip
\noindent {\bf Proof.}\quad Let $\varepsilon  >0$. Consider a small parameter $0<\delta<\varepsilon $. Let
$ (x_o,y_o)\in \R^d\times\R^d$ such that $ |x_o-y_o|<\delta$. Set
$ \eta_t = X_t(x_o)-X_t(y_o)$ and $ |\xi_t|=|\eta_t|^2$. Define
$$ \tau(x_o,y_o) = \inf\{t>0,\ \xi_t\geq \varepsilon \ \}.$$
Consider
$$ \psi_\rho(\xi) = \int_0^\xi {ds\over sr(s) + \rho}\quad\hbox{\rm and}\quad
\Phi_\rho(\xi) = e^{\psi_\rho(\xi)}.$$
As  in proof of theorem 2.2, we have for $t<\tau(x_o,y_o)$,
$$ \Phi_\rho(\xi_t)\leq \Phi_\rho(\xi_o)e^{2C\, t}.$$
Taking $\rho = |x_o-y_o|$, we get
\begin{equation}
\Phi_\rho(\xi_t)\leq e^\rho\,e^{2C\, t},\quad\hbox{\rm for  } t<\tau(x_o,y_o).
\end{equation}
Fix the point $x_o$. If $\liminf_{y_o\ra x_o}\tau(x_o,y_o) =\tau <+\infty$, we can choose 
$y_n\ra x_o$ such that 
$\lim_{n\ra +\infty}\tau(x_o,y_n) = \tau$. Applying $(17)$ for $(x_o, y_n)$ and letting 
$t\uparrow \tau(x_o,y_n)$, we get
$$ \Phi_{\rho_n}(\varepsilon )
=\Phi_{\rho_n}\bigl(\xi_{\tau(x_o,y_n)}\bigr)
\leq e^{\rho_n}\, e^{2C\,\tau(x_o,y_n)}.$$
where $\rho_n = |x_o-y_n|$. Letting $n\ra +\infty$, we have
$$ +\infty = \Phi_o(\varepsilon ) \leq e^{2C\, \tau}$$
which is absurd. Therefore
$$ \lim_{y_o\ra x_o} \tau(x_o,y_o)=+\infty,$$
which means that for any $t>0$, there exists $\delta >0$ such that for $|y_o-x_o|<\delta$,
$ \tau(x_o,y_o)>t$. In other words, 
$$ |X_t(x_o)-X_t(y_o)| \leq \varepsilon.$$

\medskip
\noindent{\bf Proposition 2.6}\quad {\it Assume that the conditions $(4)$ and $(11)$ hold. Then for $x_o\neq y_o$, we have $X_t(x_o)\neq X_t(y_o)$ for all $t\geq 0$.}

\medskip
\noindent{\bf Proof.}\quad Let $\eta_t=X_t(x_o)-X_t(y_o)$ and $\xi_t=|\eta_t|^2$. Without loss of  generality, assume that $0<\xi_o<1/2$. Let
$$\tau =\inf\bigl\{t>0,\ \xi_t\geq {3\over 4}\bigr\}.$$

By starting  from $\tau$ again , it is enough to prove that $\xi_t>0$ for $t<\tau$. Consider
$$ \psi_\rho(\xi) = \int_0^\xi {ds\over sr(s) + \rho}\quad\hbox{\rm and}\quad
\Phi_\rho(\xi) = e^{\psi_\rho(\xi)}.$$

By assumption $(10)$, for $t<\tau$, we get
$$\Bigl| {d\Phi_\rho(\xi_t)\over dt}\Bigr|
\leq 2C\Phi_\rho(\xi_t).$$
It follows that
$\Phi_\rho(\xi_t)\geq \Phi_\rho(\xi_o)-2C\int_0^t \Phi_\rho(\xi_s)ds$
or 
\begin{equation}
\Phi_\rho(\xi_t)\geq \Phi_\rho(\xi_o)e^{-2Ct}\quad \hbox{\rm for}\quad t<\tau.
\end{equation}

\noindent For $\rho>0$ small enough, 
$\Phi_\rho(\xi_o)e^{-2Ct}>1$. It follows that $\Phi_\rho(\xi_t)>1$ or $\xi_t>0$.

\medskip
Now using $(11)$ and proposition 2.6, and by the standard arguments, we obtain
\medskip

\noindent{\bf Theorem 2.7}\quad {\it Assume that the conditions $(4)$ and $(11)$ hold.
Then for any $t>0$, $x_o\ra X_t(x_o)$ defines a flow of homeomorphisms of $\R^d$.}

\vskip 10mm
\section{ Stochastic differential equations}
\bigskip

\noindent Let $\sigma: \R^d\ra \R^d\otimes\R^m$ be continuous function. Let $X(t)$ be a solution of the following It\^o stochastic differential equation:

\begin{equation}
 dX(t)=\sigma(X(t))\,dW_t + b(X(t))\, dt,\quad X(0)=x_o
\end{equation}
with the lifetime $\zeta(w)$.
\vskip 0.3cm

\noindent {\bf Theorem 3.1}\quad {\it Let $r:[1,+\infty[\ra \R_+$ be a function of ${\cal C}^1$, satisfying
\medskip
\quad (i)\qquad $\lim_{s\ra +\infty} r(s)=+\infty$,
\medskip
\quad (ii)\qquad $\int_1^\infty {ds\over sr(s)+1}=+\infty$  and

\noindent (iii)\qquad $\lim_{s\ra +\infty} {sr'(s)\over r(s)}=0$.
\vskip 2mm
Assume that for $|x|\geq 1$,
\begin{equation}
\left\{\matrix{||\sigma(x)||^2&\leq& C\, \Bigl(|x|^2\,r(|x|^2)+1\Bigr),\cr
|b(x)|&\leq& C\,  \Bigl(|x|\,r(|x|^2)+1\Bigr).}\right.
\end{equation}
Then the s.d.e $(19)$ has no explosion:
$ P(\zeta=+\infty)=1$.}

\medskip
\noindent {\bf Proof.}\quad For $0<s\leq 1$, define $r(s)=r({1\over s})$. Then there exists a constant $C>0$ such that the condition $(20)$ holds for any $x$.
Consider
$$ \psi(\xi) = \int_0^\xi {ds\over sr(s)+1}
\quad\hbox{\rm and}\quad \Phi(\xi)= e^{\psi(\xi)},\quad \xi\geq 0. $$
We have

\begin{equation}
 \Phi'(\xi)\,\bigl(\xi\,r(\xi)+1\bigr)=\Phi(\xi),
\end{equation}
\begin{equation}
 \Phi''(\xi) =\left\{\matrix{{\Phi(\xi)\bigl(1-r({1\over\xi})+{1\over\xi}r'({1\over\xi})\bigr) 
\over
(\xi r({1\over\xi})+1)^2}&\hbox{\rm if }& 0<\xi<1,\cr\cr
{\Phi(\xi)\,\bigl(1-r(\xi)-{\xi}r'({\xi})\bigr)\over (\xi r({\xi})+1)^2}&\hbox{\rm if }& \xi>1.}
\right.
\end{equation}
By conditions $(i)$ and $(iii)$, there exists $M>0$ such that
$\Phi''(\xi)\leq 0$ for $\xi\leq {1\over M}$ or $\xi\geq M$.
Remark that the function $\Phi$ is not $ {\cal C}^2$ at the point $\xi=1$. Fix a small $\delta>0$, take
$ \tilde\Phi\in {\cal C}^2(\R_+)$ such that
\begin{equation}
\tilde\Phi\geq \Phi,\quad \tilde\Phi(\xi)=\Phi(\xi)\quad
\hbox{\rm for  }\xi\not\in [1-\delta, 1+\delta].
\end{equation}
Denote
$$ K_1=\sup_{\xi\in [1-\delta, 1+\delta]}
\Bigl( |\tilde\Phi'(\xi)|+|\tilde\Phi''(\xi)|\Bigr),
\quad K_2 = \sup_{\xi\in [1-\delta, 1+\delta]}
\Bigl( \xi\,r({\xi})\Bigr).$$
Then  
\begin{equation}
 |\tilde\Phi'(\xi)|\leq {K_1(K_2+1)\over \Phi(1-\delta)}
\cdot {\Phi(\xi)\over \xi\, r({\xi})+1 },\quad \xi\in [1-\delta, 1+\delta],
\end{equation}
and
\begin{equation}
|\tilde\Phi''(\xi)|\leq {K_1(K_2+1)^2\over \Phi(1-\delta)}
\cdot {\Phi(\xi)\over (\xi\, r({\xi})+1)^2 }, \quad \xi\in [1-\delta, 1+\delta].
\end{equation}
Let  $ \xi_t(w)=|X_t(w)|^2$. We have

\begin{eqnarray}
d\xi_t&=&2\bigl<X_t, \sigma(X(t))\,dW_t\bigr>
+ 2\bigl<X_t, b(X(t))\bigr>\,dt\nonumber \\
&& \hskip 10mm +||\sigma(X(t))||^2\, dt,
\end{eqnarray}
and the stochastic contraction $d\xi_t\cdot d\xi_t$ is given by
\begin{equation}
 d\xi_t\cdot d\xi_t
=4|\sigma^*(X_t)X_t|^2\,dt
\end{equation}
where $\sigma^*$ denotes the transpose matrix of $\sigma$.

\noindent Define
$$ \tau_R =\inf{\bigl\{t>0, \ \xi_t\geq R\,\bigr\}},\quad R>0.$$
Then $\tau_R\uparrow\zeta$ as $ R\uparrow +\infty$. Let
$$I_w=\bigl\{t>0,\ \xi_t(w)\in [{1\over M}, M]\bigr\}.$$
Taking $M$ big enough such that $[1-\delta, 1+\delta]\subset [{1\over M}, M]$. By $(22)$,
\begin{equation}
 \tilde\Phi''(\xi_t)=\Phi''(\xi_t)\leq 0\quad\hbox{\rm for  } t\not\in I_w.
\end{equation}
Combining $(22)$ and $(25)$, there exists a constant $C_1$ such that
\begin{equation}
 \Bigl|\tilde\Phi''(\xi_t)\Bigr|
\leq {C_1\Phi(\xi_t)\over (\xi_t\, r({\xi_t})+1)^2 },\quad t\in I_w.
\end{equation}
By $(21)$ and $(24)$, for some constant $C_2>0$, we have
\begin{equation}
\Bigl|\tilde\Phi'(\xi_t)\Bigr|
\leq {C_2\,\Phi(\xi_t)\over \xi_t\, r({\xi_t})+1},\quad t>0.
\end{equation}
Now by It\^o formula and according to $(26), (27)$, we have
\begin{eqnarray}
\tilde\Phi(\xi_{t\wedge\tau_R})
&=&\Phi(\xi_o) + 2\int_0^{t\wedge\tau_R} \tilde\Phi'(\xi_s)
\bigl<X_s, \sigma(X(s))dW_s\bigr>\nonumber\\
&+& 2\int_0^{t\wedge\tau_R} \tilde\Phi'(\xi_s)\,
\bigl<X_s, b(X(s))\bigr>\,ds\nonumber\\
&+&\int_0^{t\wedge\tau_R}\tilde\Phi'(\xi_s)
||\sigma(X(s))||^2\,ds\nonumber\\
&+ &2 \int_0^{t\wedge\tau_R}\tilde\Phi''(\xi_s)\,
|\sigma^*(X(s))X_s|^2\,ds.
\end{eqnarray}
By $(28)$ and $(29)$,
\begin{equation}
 \int_0^{t\wedge\tau_R}\tilde\Phi''(\xi_s)\,
|\sigma^*(X(s))X_s|^2\,ds
\leq \int_0^{t\wedge\tau_R}{\bf 1}_{I_w}(s){C_1\,\Phi(\xi_s)\over (\xi_s\,r({\xi_s})+1)^2 }\,|\sigma^*(X(s))X_s|^2\,ds.
\end{equation}
By $(20)$, 
\begin{equation}
{|\sigma^*(X(s))X_s|^2\over (\xi_s\, r({\xi_s})+1)^2 }
\leq C_3\ {\xi_s\, (\xi_s\, r(\xi_s)+1)\over (\xi_s\, r(\xi_s)+1)^2 }
\end{equation}
which is dominated by a constant $C_3$. According to $(32)$, we get

\begin{equation}
 \int_0^{t\wedge\tau_R}\tilde\Phi''(\xi_s)\,|\sigma^*(X(s))X_s|^2\,ds
\leq C_3 \int_0^{t\wedge\tau_R} \Phi(\xi_s)\, ds.
\end{equation}
In the same way, for some constant $C_4>0$, we have

\begin{equation}
 {|\bigl<X_s, b(X(s))\bigr>| + ||\sigma(X(s))||^2\over
\xi_s\, r({\xi_s}) + 1}\leq C_4, \quad s>0.
\end{equation}
Now using $(31)$ and according to $(30)$, $(35)$ and $(34)$, we get

$$ E\Bigl(\Phi(\xi_{t\wedge\tau_R})\Bigr)
\leq E\Bigl(\tilde\Phi(\xi_{t\wedge\tau_R})\Bigr)
\leq \Phi(\xi_o)+ C_5\,\int_0^t E\bigl(\Phi(\xi_{s\wedge\tau_R})\bigr)\, ds,$$
which implies that

$$ E\Bigl(\Phi(\xi_{t\wedge\tau_R})\Bigr)
\leq \Phi(\xi_o)e^{C_5\, t}.$$
Letting $R\ra +\infty$, by Fatou lemma, we get
\begin{equation}
 E\Bigl(\Phi(\xi_{t\wedge\zeta})\Bigr)
\leq \Phi(\xi_o)e^{C_5\, t}.
\end{equation}
Now if $ P(\zeta <+\infty) >0$, then for some $T>0$, $P(\zeta \leq T)>0$. Taking $t=T$ in $(36)$, we get
$$ E\Bigl( {\bf 1}_{(\zeta\leq T)}\Phi(\xi_{\zeta})\Bigr)
\leq \Phi(\xi_o)e^{C_5\, t}$$
which is impossible, because of $\Phi(\xi_\zeta)=+\infty$.

\medskip
\noindent {\bf Theorem 3.2}\quad {\it  Let $r:\ ]0,1[\ra \R_+$ be ${\cal C}^1$-function satisfying the conditions
\medskip

\quad (i)\quad $\lim_{\xi\ra 0} r(\xi)=+\infty$,

\smallskip
\quad (ii)\quad $ \lim_{\xi\ra 0}{\xi\, r'(\xi)\over r(\xi)}=0$,

\smallskip
\quad (iii)\quad $\int_0^a {ds\over sr(s)} =+\infty$, for any $a>0$.
\medskip

Assume that for $|x-y|<1$,
\begin{equation}
 \left\{\matrix{||\sigma(x)-\sigma(y)||^2 &\leq & C\, |x-y|^2\,r(|x-y|^2),\cr
        |b(x)-b(y)| &\leq& C\,  |x-y|\,r(|x-y|^2).}\right.
\end{equation}
Then the s.d.e. $(19)$ has the pathwise uniqueness.}

\medskip

\noindent {\bf Proof.}\quad Let $ \eta_t(w)=X(t)-Y(t)$ and 
$ \xi_t(w) = |\eta_t(w)|^2$. 
Let $\rho>0$. Define the function $\psi_\rho: [0,1]\ra \R$ by
\begin{equation}
 \psi_\rho(\xi) =\int_0^\xi {ds\over sr(s)+\rho}.
\end{equation}
It is clear that for any $0<\xi<1$, 
$$ \psi_\rho(\xi)\uparrow \psi_0(\xi)=\int_0^\xi {ds\over sr(s)}=+\infty,
\quad\hbox{\rm as  } \rho\downarrow 0.$$
Define
\begin{equation}
 \Phi_\rho(\xi) =e^{\psi_\rho(\xi)}.
\end{equation}
Then we have
\begin{equation}
 \Phi_\rho'(\xi)\, (\xi\,r(\xi)+\rho) =\Phi_\rho(\xi),
\end{equation}
and
\begin{equation}
 \Phi_\rho''(\xi) ={\Phi_\rho(\xi)(1-\xi\,r'(\xi)-r(\xi))\over
(\xi\,r(\xi)+\rho)^2}.
\end{equation}
By assumption $(i), (ii)$ on $r$, there exists $\delta >0$ such that 
$\Phi_\rho''(\xi)\leq 0$ for $0<\xi<\delta$.

\noindent Let 
$$\tau =\inf{\bigl\{\, t>0,\quad \xi_t\geq \delta\,\bigr\}}.$$
By It\^o formula, we have
\begin{eqnarray}
\Phi_\rho(\xi_{t\wedge\tau})
&=&1 + 2\int_0^{t\wedge\tau} \Phi_\rho'(\xi_s)
\bigl<\eta_s, \bigl(\sigma(X(s))-\sigma(Y(s))\bigr)dW_s\bigr>\nonumber\\
&+& 2\int_0^{t\wedge\tau} \Phi_\rho'(\xi_s)\,
\bigl<\eta_s, b(X(s))-b(Y(s))\bigr>\,ds\nonumber\\
&+&\int_0^{t\wedge\tau}\Phi_\rho'(\xi_s)
||\sigma(X(s))-\sigma(Y(s))||^2\,ds\nonumber\\
&+ & 2 \int_0^{t\wedge\tau}\Phi_\rho''(\xi_s)\,
|\bigl(\sigma^*(X(s))-\sigma^*(Y(s))\bigr)\eta_s|^2\,ds.
\end{eqnarray}
Applying the hypothesis $(37)$, we get

$$ E\Bigl(\Phi_\rho(\xi_{t\wedge\tau})\Bigr)
\leq 1 + 2C\, E\Bigl(\int_0^{t\wedge\tau} \Phi_\rho'(\xi_{s})\,\xi_{s}\,r(\xi_s)\ ds\Bigr)$$
which is smaller by $(41)$ than
$$ 1 + 2C \int_0^t E\Bigl(\Phi_\rho(\xi_{s\wedge\tau})\Bigr)\, ds.$$
Now by Gronwall lemma, we get
$ E\Bigl(\Phi_\rho(\xi_{t\wedge\tau})\Bigr)
\leq e^{2ct}$ or
\begin{equation}
 E\Bigl(e^{\psi_\rho(\xi_{t\wedge\tau})}\Bigr)\leq e^{2Ct}.
\end{equation}
Letting $\rho\downarrow 0$ in $(43)$, 
$$ E\Bigl(e^{\psi_0(\xi_{t\wedge\tau})}\Bigr)\leq e^{2Ct}$$
which implies that for any $t$ given,  
\begin{equation}
 \xi_{t\wedge\tau}=0\quad\hbox{\rm almost surely.}
\end{equation}
If
$ P(\tau<+\infty)>0$, then for some $T>0$ big enough
$P(\tau\leq T)>0$. By $(44)$, almost surely for all $ t\in Q\cap [0,T]$, 
$\xi_{t\wedge\tau}=0$. It follows that on $\bigl\{\tau\leq T\bigr\}$,
$$\xi_\tau = 0$$
which is absurd with the definition of $\tau$. Therefore $\tau=+\infty$ almost surely and for 
any $t$ given, $\xi_t=0$ almost surely. Now by continuity of samples, the two solutions are indistinguishable.

\medskip
\noindent {\bf Remark 3.3}\quad In the case of $d=m=1$, finer results about pathwise uniqueness have been established. Namely $\sigma$ was allowed to be H\"older of exponent $\geq 1/2$ (see [RY, Ch. IX-3], [IW, p.168]).

\medskip
\noindent {\bf Theorem 3.4}\quad {\it Assume that the coefficients $\sigma$ and $b$ satisfy the assumption $(20)$ and $(37)$. Let $ X(t,x_o)$ be the solution of the s.d.e. $(19)$ with initial value $x_o$. Then for any $\varepsilon>0$, we have}
\begin{equation}
 \lim_{y_o\ra x_o} P( \sup_{0\leq s\leq t}|X(s,y_o)-X(s,x_o)|>\varepsilon )=0.
\end{equation}

\medskip
\noindent {\bf Proof.}\quad Fix $x_o$. Let $\delta$ be the parameter given in proof of theorem 3.2, consider $|y_o-x_o|<\varepsilon<\delta$. Let 
$\xi_t(w)=|X(t,y_o)-X(t,x_o)|^2$. Define
$$ \tau_w(x_o,y_o)=\inf\{\, t>0,\ \xi_t>\varepsilon^2\,\}.$$
The same arguments as above yields to
$$ E\Bigl(\Phi_\rho(\xi_{t\wedge\tau(x_o,y_o)})\Bigr) \leq \Phi_\rho(\xi_o)\, e^{2Ct}.$$
Taking $\rho=|x_o-y_o|$, we have 
$ E\Bigl(\Phi_\rho(\xi_{t\wedge\tau(x_o,y_o)})\Bigr)
\leq e^\rho\, e^{2Ct}$. Hence
$$ P(\tau(x_o,y_o)<t)\Phi_\rho(\varepsilon )
\leq E\Bigl(\Phi_\rho(\xi_{t\wedge\tau(x_o,y_o)})\Bigr)\leq e^\rho\, e^{2Ct}.$$
It follows
$$ P(\sup_{0\leq s\leq t}|X(s,y_o)-X(s,x_o)|>\varepsilon )
=P(\tau(x_o,y_o)<t)\leq e^{-\psi_\rho(\varepsilon )}\,e^\rho\, e^{2Ct}\ra 0$$
as $\rho=|y_o-x_o|\ra 0$.

\vskip 10mm
\noindent{\bf Corollary 3.5} {\it The diffusion process $(X(t,x))$ given by the solution of the s.d.e. is Feller, i.e., the associated semigroup $(T_t, t\geq 0)$ maps $C_b(R^d)$ into $C_b(R^d)$.}
\vskip 0.3cm
\noindent{\bf Proof}. It is a direct consequence of Theorem 3.4 and the definition
$$T_tf(x)=E[f(X(t,x))], \quad \quad f\in  C_b(R^d)$$

\medskip
\noindent {\bf Remark 3.6}\quad We have a difficulty here to apply the Kolmogoroff modification theorem to obtain a version $\tilde X(t,x_o)$ such that $x_o\ra \tilde X(t,x_o)$ is continuous. However the situation for the case of $S^1$ is well handled (see [F], [M]).

\section{Statement of large deviation principle}
\bigskip

\noindent  Let $ \sigma:\ R^d\rightarrow  R^d\otimes R^m$ and 
$ b:\ R^d\rightarrow  R^d$ be continuous  functions. Consider the following It\^o s.d.e:

\begin{equation}
 dX^{\varepsilon}(t)=\varepsilon^{\frac{1}{2}}\sigma(X^{\varepsilon}(t))\,dW_t + b(X^{\varepsilon}(t))\, dt,\quad X_o^{\varepsilon}(w)=x_o
\end{equation}
where $t\rightarrow  W_t$ is a $R^m$-valued standard Brownian motion.  In order to be more explicit, we shall work under the following assumptions,

$$ \left\{\matrix{||\sigma(x)-\sigma(y)||^2&\leq& C\, |x-y|^2\,\log{{ 1\over
 |x-y|}},
&\quad\hbox{\rm for  } |x-y|<1,\cr
|b(x)-b(y)|&\leq& C\, |x-y|\,\log{{ 1\over  |x-y|}},&\quad\hbox{\rm for  }
 |x-y|<1}\right.
\leqno(H1)$$

$$ \left\{\matrix{||\sigma(x)||^2 &\leq & C\,( |x|^2\,\log{|x|}+1),\cr
|b(x)| &\leq & C\,(|x|\,\log{|x|}+1)}\right.
\leqno(H2)$$
where $|\cdot|$ denotes the Euclidean distance in $R^d$ and 
$ ||\sigma||^2=\sum_{i=1}^d\sum_{j=1}^m \sigma_{ij}^2$.  The rest of this paper  is to establish a large deviation principle for solutions of above s. d. e.
\vskip 0.3cm
\noindent Let  $ C_x\big([0,1],R^m\big)$ be the space of continuous functions from $[0,1]$ into $R^m$ with initial value $x$. If $g\in C_0\big([0,1],R^m\big)$  is absolutely continuous, set $e(g)=\int_0^1 |\dot{g}(t)|^2 dt$. Let $F(g)$ be the solution to the differential equation
\begin{eqnarray}
F(g)(t) &=& x_0+\int_0^t b\big( F(g)(s)\big)ds \nonumber \\
	& & +\int_0^t \sigma \big( F(g)(s)\big)\dot{g}(s)ds, \quad 0 < t < \infty
\end{eqnarray}
the uniqueness and non explosion being obtained as in section 2 for the differential equation $(2)$. See also lemma 5.3 below.
\vskip 0.3cm

\noindent {\bf Theorem 4.1}   \label{thm:threeone} 
	{\it Let $\mu_{\varepsilon}$ be the law of $X^{\varepsilon}(\cdot)$ on 
$C_{x_0}\big([0,1], R^d\big)$. Assume $ (H1) $ and $(H2)$. Then $\{\mu_{\varepsilon},\varepsilon>0 \}$ 
satisfies a large deviation principle with the following  good rate function 
\begin{equation} \label{eq:eight}  
  I(f)=\inf\bigg\{\frac{1}{2}e(g)^2; F(g)=f\bigg\}, \quad \quad 
	f \in C_{x_0}([0,1], R^d)
\end{equation}
i.e., 
\begin{enumerate}
\item[(i)] for any closed subset $C\subset  C_{x_0}([0,1], R^d)$,
\begin{equation} \label{eq:nine}  
  \limsup_{\varepsilon \rightarrow 0} \varepsilon \log\mu_{\varepsilon}(C)\leq -\inf_{f\in C}I(f)
\end{equation}
\item[(ii)] for any open  subset $G\subset  C_{x_0}([0,1], R^d)$,
\begin{equation}  \label{eq:ten}  
  \liminf_{\varepsilon \rightarrow 0} \varepsilon \log\mu_{\varepsilon}(G)\geq  -\inf_{f\in G}I(f)
\end{equation}
\end{enumerate}}
\noindent The proof of the theorem will be given in Section 5 and 6. 
 
\section{Large deviations when $\sigma$, $b$ are bounded }

The theory of large deviations for diffusion processes under Lipschitzian coefficients is well established (see [A], [S]). Some new developments in infinite dimensional situations are discussed in [FZ1,2], [Z1,2]. The main task of this work is to handle the non Lipschitzian feature. For  $n\geq 1$, let  $X_n^{\varepsilon}(\cdot)$ to be the solution to
\begin{eqnarray} \label{eq:eleven}  
X_n^{\varepsilon}(t) &= &x_0
	     +\int_0^t b(X_n^{\varepsilon}(\frac{[ns]}{n}))ds \nonumber \\
	& &  + \varepsilon^{\frac{1}{2}} \int_0^t \sigma \left(X_n^{\varepsilon}\bigg(\frac{[ns]}{n}\bigg)\right)dW_s	
\end{eqnarray}

\noindent We need the following lemma from Stroock [S,P.81].
\vskip 0.3cm

\noindent {\bf Lemma 5.1}.  \label{lem:threeone}
	{\it Let $\alpha (\cdot )$ and $\beta (\cdot )$ be $({\mathcal F}_t)_{t \geq 0}$ -progressively measurable with values in $R^d\otimes R^m$ and $R^d$ respectively. Assume $\vert\vert \alpha (\cdot )\vert \vert \leq A$ and $\vert \beta (\cdot )\vert  \leq B$, and set $\xi (t)=\int_0^t \alpha (s)dW_s +\int_0^t \beta (s) ds$. Then for $T>0$ and $R>0$ satisfying $d^{\frac{1}{2}}BT<R$:
\begin{equation} \label{eq:twelve}  
  P\big(\sup_{0\leq t\leq T}|\xi (t)|\geq R\big)\leq 2d  
		\exp\big(-(R-d^{\frac{1}{2}}BT)^2/2A^2dT\big)
\end{equation}}

\noindent {\bf Proposition 5.2}. \label{prop:threeone}
	{\it In addition to (H.1), we also assume that $b,\sigma$ are bounded. For any $\delta>0$, it holds that
\begin{equation} \label{eq:thirteen}  
  \lim_{n\rightarrow \infty}\limsup_{\varepsilon \rightarrow 0} \varepsilon 
	\log P\big(\sup_{0\leq t\leq 1}\vert X^{\varepsilon}(t)
	-X_n^{\varepsilon}(t)\vert > \delta \big)=-\infty
\end{equation}}

\noindent {\bf Proof}. We may and will assume $\delta <e^{-1}<1$. Let  
$Y_n^{\varepsilon}(t):= X^{\varepsilon}(t)-X_n^{\varepsilon}(t)$ and $\xi_n^{\varepsilon}(t)=|Y_n^{\varepsilon}(t)|^2$.  We have 
\begin{eqnarray} \label{eq:fourteen}
  Y_n^{\varepsilon}(t) 
	&= &\int_0^t\left[ b\big(X^{\varepsilon}(s)\big)
	     - b\big(X_n^{\varepsilon}(\frac{[ns]}{n})\big)\right]ds\nonumber   \\
	& &  + \varepsilon^{\frac{1}{2}}\int_0^t\left[ \sigma\big(X^{\varepsilon}(s)\big)
	     - \sigma\big(X_n^{\varepsilon}(\frac{[ns]}{n})\big)\right] dW_s 
\end{eqnarray}
and 
\begin{eqnarray}
d\xi_n^{\varepsilon}(t)&=& 2 \varepsilon^{\frac{1}{2}}\big <Y_n^{\varepsilon}(t),\left( \sigma\big(X^{\varepsilon}(t)\big)
	     - \sigma\big(X_n^{\varepsilon}(\frac{[nt]}{n})\big)\right) dW_t\big >\nonumber \\
 & & +2\big <Y_n^{\varepsilon}(t),\left( b\big(X^{\varepsilon}(t)\big)
	     - b\big(X_n^{\varepsilon}(\frac{[nt]}{n})\big)\right)\big >dt\nonumber \\
&& +\varepsilon ||\sigma\big(X^{\varepsilon}(t)\big)
	     - \sigma\big(X_n^{\varepsilon}(\frac{[nt]}{n})\big)||^2dt
\end{eqnarray}

 \noindent The stochastic contraction $d\xi_n^{\varepsilon} \cdot d\xi_n^{\varepsilon}$ is given by
\begin{equation}
 d\xi_n^{\varepsilon} \cdot d\xi_n^{\varepsilon}
=4\varepsilon |\left( \sigma\big(X^{\varepsilon}(t)\big)
	     - \sigma\big(X_n^{\varepsilon}(\frac{[nt]}{n})\big)\right)^*Y_n^{\varepsilon}(t)|^2\,dt
\end{equation}
where $\sigma^*$ denotes the transpose  of $\sigma$.
Let $\rho>0$. Define the function $\psi_\rho: [0,1]\rightarrow  R$ by
\begin{equation}
 \psi_\rho(\xi) =\int_0^\xi {ds\over s\log{{1\over s}}+\rho}.
\end{equation}
It is clear that for any $0<\xi<1$, 
$$ \psi_\rho(\xi)\uparrow \psi_0(\xi)=\int_0^\xi {ds\over s\log{{1\over s}}}=+\infty,
\quad\hbox{\rm as  } \rho\downarrow 0.$$
Define for $\lambda>0$
\begin{equation}
 \Phi_{\rho,\lambda}(\xi) =e^{\lambda \psi_\rho(\xi)}.
\end{equation}
Then we have
\begin{equation}
 \Phi_{\rho, \lambda}'(\xi)\, (\xi\log{{1\over\xi}}+\rho) =\lambda \Phi_{\rho, \lambda}(\xi),
\end{equation}
and
$$
 \Phi_{\rho,\lambda}''(\xi) =\lambda^2 \Phi_{\rho, \lambda}(\xi){1\over
(\xi\log{{1\over\xi}}+\rho)^2}+\lambda \Phi_{\rho, \lambda}(\xi){(1-\log{{1\over\xi}})\over
(\xi\log{{1\over\xi}}+\rho)^2}
$$
\begin{equation}
 \leq\lambda^2 \Phi_{\rho, \lambda}(\xi){1\over
(\xi\log{{1\over\xi}}+\rho)^2}  \quad\hbox{\rm if  }\xi\leq e^{-1}.
\end{equation}
Now, choose  a positive constant $\delta_1<e^{-1}$ satisfying $\delta_1\log{\frac{1}{\delta_1}}<\rho$.    Define $\tau_{n}^{\varepsilon}=\inf\{ t\geq 0; \vert X_n^{\varepsilon}(t)-X_n^{\varepsilon}(\frac{[nt]}{n})\vert \geq \delta_1\}$, and set $\xi_{n,\delta_1}^{\varepsilon }(t)=\xi_n^{\varepsilon}(t\wedge  \tau_{n}^{\varepsilon}),t\geq 0$. Putting  $T_{n}^{\varepsilon}=\inf\{t\geq 0, \vert \xi_{n,\delta_1}^{\varepsilon }(t)\vert \geq \delta^2 \}$, we have,
\begin{eqnarray} \label{eq:fifteen}
	P\big(\sup_{0\leq t\leq 1}\vert Y_n^{\varepsilon}(t)\vert > \delta \big)
	&=& P\big(\sup_{0\leq t\leq 1}\vert Y_n^{\varepsilon}(t)\vert 
		> \delta ,\tau_{n}^{\varepsilon}\leq 1 \big)\nonumber	\\
 	& &  + P\big(\sup_{0\leq t\leq 1}\vert Y_n^{\varepsilon}(t)\vert 
		> \delta ,\tau_{n}^{\varepsilon}> 1 \big)\nonumber	\\
 	&\leq & P(\tau_{n}^{\varepsilon}\leq 1 )+P(T_{n}^{\varepsilon} \leq 1)
\end{eqnarray}
\noindent Observe,
\begin{equation} \label{eq:sixteen}
  P(\tau_{n}^{\varepsilon}\leq 1 )\leq \sum_{k=1}^n P\bigg(\sup_{\frac{k-1}{n}\leq t 
	\leq \frac{k}{n}}\vert X_n^{\varepsilon}(t)-X_n^{\varepsilon}(\frac{k-1}{n})\vert 
	\geq \delta_1 \bigg).
\end{equation}
\noindent Using Lemma 5.1 and the boundness of $\sigma, b$, there exists a constant $c_{\delta_1} >0$ such that 
\begin{equation} \label{eq:seventeen}
  P(\tau_{n}^{\varepsilon}\leq 1 )\leq n\exp(-nc_{\delta_1}/\varepsilon ).
\end{equation}
\noindent Hence,
\begin{equation} \label{eq:eightteen}
  \lim_{n\rightarrow \infty}\limsup_{\varepsilon \rightarrow 0} 
	\varepsilon \log P(\tau_{n,\rho}^{\varepsilon}\leq 1 )=-\infty
\end{equation}

\noindent For notational simplicity, write $T$ for $T_{n}^{\varepsilon}$ and $\tau$ for $\tau_{n}^{\varepsilon}$.  By Ito's formula,

\begin{eqnarray}
&\Phi_{\rho,\lambda}(\xi_{n,\delta_1}^{\varepsilon }(t\wedge T))\nonumber\\
=&1 + 2\varepsilon^{1\over 2}\int_0^{t\wedge T\wedge \tau } \Phi_{\rho,\lambda}'(\xi_{n,\delta_1}^{\varepsilon }(s))
\bigg <Y_n^{\varepsilon}(s), \left( \sigma (X^{\varepsilon}(s))
	     - \sigma(X_n^{\varepsilon}(\frac{[ns]}{n}))\right) dW_s\bigg >\nonumber \\
+ &2\int_0^{t\wedge T\wedge \tau } \Phi_{\rho,\lambda}'(\xi_{n,\delta_1}^{\varepsilon }(s))\,
\bigg <Y_n^{\varepsilon }(s),  b(X^{\varepsilon }(s))
	     - b(X_n^{\varepsilon }(\frac{[ns]}{n})) \bigg > ds\nonumber \\
+&\varepsilon \int_0^{t\wedge T\wedge \tau }\Phi_{\rho,\lambda}'(\xi_{n,\delta_1}^{\varepsilon }(s))
|| \sigma(X^{\varepsilon}(s))
	     - \sigma(X_n^{\varepsilon}(\frac{[ns]}{n}))||^2 ds\nonumber\\
+& 2 \varepsilon \int_0^{t\wedge T\wedge \tau }\Phi_{\rho,\lambda }^{\prime\prime }(\xi_{n,\delta_1 }^{\varepsilon }(s))\,
|\left( \sigma\big(X^{\varepsilon }(s)\big)
	     - \sigma\big(X_n^{\varepsilon}(\frac{[ns]}{n})\big)\right)^*Y_n^{\varepsilon}(s)|^2\,ds.
\end{eqnarray}
Note that for $s\leq T\wedge \tau$, $|Y_n^{\varepsilon}(s)|\leq \delta_1 \leq e^{-1}<1$ and $|X_n^{\varepsilon}(s)-X_n^{\varepsilon}(\frac{[ns]}{n})|\leq \delta_1<e^{-1}$. Therefore, for $s\leq T\wedge\tau$,
$$|\bigg <Y_n^{\varepsilon}(s),  b\big(X^{\varepsilon}(s)\big)
	     - b\big(X_n^{\varepsilon}(\frac{[ns]}{n})\big) \bigg >|$$
$$=|\bigg <Y_n^{\varepsilon}(s),  b\big(X^{\varepsilon}(s)\big)
	     - b\big(X_n^{\varepsilon}(s)\big) \bigg >|+
|\bigg <Y_n^{\varepsilon}(s),  b\big(X_n^{\varepsilon}(s)\big)
	     - b\big(X_n^{\varepsilon}(\frac{[ns]}{n})\big) \bigg >|$$
$$\leq C|Y_n^{\varepsilon}(s)|^2 \log(\frac{1}{|Y_n^{\varepsilon}(s)|})+C |X_n^{\varepsilon}(s)-X_n^{\varepsilon}(\frac{[ns]}{n})|\log(\frac{1}{|X_n^{\varepsilon}(s)-X_n^{\varepsilon}(\frac{[ns]}{n})|})$$
\begin{equation}
\leq \frac{1}{2}C \xi_n^{\varepsilon}(s)\log(\frac{1}{\xi_n^{\varepsilon}(s)})+ \delta_1\log(\frac{1}{\delta_1})\leq C\big ( \xi_n^{\varepsilon}(s)\log(\frac{1}{\xi_n^{\varepsilon}(s)})+\rho\big ),
\end{equation}
where we have used the fact that the function  $x\log(\frac{1}{x})$ is increasing on $[0, e^{-1}]$.
Furthermore,for $s\leq T\wedge\tau$,
$$||\sigma\big(X^{\varepsilon}(s)\big)
	     - \sigma\big(X_n^{\varepsilon}(\frac{[ns]}{n})\big) ||^2$$
$$\leq 2||\sigma\big(X^{\varepsilon}(s)\big)
	     - \sigma\big(X_n^{\varepsilon}(s))\big) ||^2+2||\sigma\big(X_n^{\varepsilon}(s)\big)
	     - \sigma\big(X_n^{\varepsilon}(\frac{[ns]}{n})\big)||^2$$
$$
\leq 2C \xi_n^{\varepsilon}(s)\log(\frac{1}{\xi_n^{\varepsilon}(s)})+2C |X_n^{\varepsilon}(s)-X_n^{\varepsilon}(\frac{[ns]}{n})|^2\log(\frac{1}{|X_n^{\varepsilon}(s)-X_n^{\varepsilon}(\frac{[ns]}{n})|})
$$
\begin{equation}
\leq C\big ( \xi_n^{\varepsilon}(s)\log(\frac{1}{\xi_n^{\varepsilon}(s)})+\rho\big )
\end{equation}
Similarly, 
$$|\left( \sigma\big(X^{\varepsilon}(s)\big)
	     - \sigma\big(X_n^{\varepsilon}(\frac{[ns]}{n})\big)\right)^*Y_n^{\varepsilon}(s)|^2$$
$$\leq |Y_n^{\varepsilon}(s)|^2||\left( \sigma\big(X^{\varepsilon}(s)\big)
	     - \sigma\big(X_n^{\varepsilon}(\frac{[ns]}{n})\big)\right)^* ||^2$$
\begin{equation}
\leq C\xi_n^{\varepsilon}(s)\big ( \xi_n^{\varepsilon}(s)\log(\frac{1}{\xi_n^{\varepsilon}(s)})+\rho\big )
\leq C\big ( \xi_n^{\varepsilon}(s)\log(\frac{1}{\xi_n^{\varepsilon}(s)})+\rho\big )^2
\end{equation}
Taking these inequalities into account, it follows from (65) that
\begin{eqnarray}
& E[\Phi_{\rho,\lambda}(\xi_{n,\delta_1}^{\varepsilon }(t\wedge T))]\nonumber\\
\leq& C E\bigg [ \int_0^{t\wedge T} \Phi_{\rho,\lambda}^{\prime}(\xi_{n,\delta_1}^{\varepsilon }(s))(\xi_{n,\delta_1}^{\varepsilon }(s)\log(\frac{1}{\xi_{n,\delta_1}^{\varepsilon }(s)})+\rho)ds\bigg ]\nonumber \\
+& C\varepsilon  E\bigg [ \int_0^{t\wedge T} \Phi_{\rho,\lambda}^{\prime\prime }(\xi_{n,\delta_1}^{\varepsilon }(s))(\xi_{n,\delta_1}^{\varepsilon }(s)\log(\frac{1}{\xi_{n,\delta_1}^{\varepsilon }(s)})+\rho)^2ds\bigg ]
\end{eqnarray}
which is smaller by (59) and (60) than
$$ 1+C\varepsilon (\lambda^2 +\lambda) \int_0^t  E[\Phi_{\rho,\lambda}(\xi_{n,\delta_1}^{\varepsilon }(s\wedge T))]ds$$
By Gronwall's inequality, we obtain that
\begin{equation}
E[\Phi_{\rho,\lambda}(\xi_{n,\delta_1}^{\varepsilon }(1\wedge T))]\leq e^{C(\varepsilon \lambda^2+\lambda)}
\end{equation}
On the other hand,
\begin{equation}
E[\Phi_{\rho,\lambda}(\xi_{n,\delta_1}^{\varepsilon }(1\wedge T))]\geq E[\Phi_{\rho,\lambda}(\xi_{n,\delta_1}^{\varepsilon }(T)), T \leq 1]=e^{\lambda \psi_{\rho} (\delta^2)}P(T\leq 1)
\end{equation}
Combining (70) with (71), we have
$$P(T\leq 1 )\leq e^{-\lambda \psi_\rho(\delta^2)}e^{C(\varepsilon \lambda^2+\lambda)}$$
Taking $\lambda=\frac{1}{\varepsilon}$ it follows that
$$\limsup_{n\rightarrow \infty}\limsup_{\varepsilon \rightarrow 0}\varepsilon logP(T\leq 1 )\leq -\psi_{\rho} (\delta^2)+2C$$
This together with (64) implies that 
\begin{equation}
\limsup_{n\rightarrow \infty}\limsup_{\varepsilon \rightarrow 0}\varepsilon logP(\sup_{0\leq t\leq 1}|Y_n^{\varepsilon}(t)|>\delta )\leq -\psi_{\rho} (\delta^2)+2C
\end{equation}
Sending $\rho$ to $0$ completes the proof.
  
\vskip 0.3cm

\noindent For $n\geq 1$, define a map $F_n(\cdot ): C_0([0,1], R^m)\rightarrow C_{x_0}([0, 1],R^d)$ by
\begin{eqnarray} \label{eq:twentynine}
F_n(\omega)(0)&=&x_0\nonumber\\
  F_n(\omega)(t) &=& F_n(\omega)\bigg(\frac{k}{n}\bigg)
		+ b\bigg( F_n(\omega)\bigg(\frac{k}{n}\bigg)\bigg)\left(t
		-\frac{k}{n}\right)\nonumber \\
	& & + \sigma \bigg( F_n(\omega)\bigg(\frac{k}{n}\bigg)\bigg)\left(\omega (t)
		-\omega \bigg(\frac{k}{n}\bigg)\right)
\end{eqnarray}
for $\frac{k}{n}\leq t\leq \frac{k+1}{n}$. 
It is easy to see that $F_n$ is a  continuous map from $ C_0([0,1], R^m)$ to $C_{x_0}([0, 1],R^d)$.


\bigskip

\noindent {\bf Lemma 5.3}.  \label{lem:threetwo}
	$\lim_{n\rightarrow \infty}\sup_{\{g;e(g)\leq \alpha\}}
\sup_{0\leq t\leq 1}|F_n(g)(t)-F(g)(t)|=0$.

\vskip 0.3cm

\noindent {\bf Proof }. Note that for $g$ with $e(g)\leq \alpha$,
\begin{eqnarray} \label{eq:thirty}
  F_n(g)(t)& = &x_0+\int_{0}^t b\big(F_n(g)(\frac{[ns]}{n})\big)ds\nonumber \\
	& + &\int_0^t \sigma \bigg(F_n(g)\bigg(\frac{[ns]}{n}\bigg)\bigg)\dot{g}(s)ds
\end{eqnarray}
\noindent Thus,
$$
    F_n(g)(t)-F(g)(t)
$$
$$
= \int_{0}^t\big[ b(F_n(g)\bigg(\frac{[ns]}{n}\bigg)\big)-b(F(g)(s))\big]ds 
$$
\begin{equation}
+ \int_0^t\big[ \sigma (F_n(g)\bigg(\frac{[ns]}{n}\bigg)\big)-\sigma (F(g)(s))\big] \dot{g}(s)ds 
\end{equation}
\noindent Since $b, \sigma$ are bounded, we have for $t\leq 1$
$$
    \big|F_n(g)(t)-F_n(g)\bigg(\frac{[nt]}{n}\bigg)\big|
		\leq \int_{\frac{[nt]}{n}}^t| b(F_n(g)
		(\frac{[ns]}{n}\big ))|ds 
$$
$$
+\int_{\frac{[nt]}{n}}^t \big||\sigma (F_n(g)
		\bigg(\frac{[ns]}{n}\bigg) )\big|||\dot{g}|(s)ds
$$
\begin{equation}
 \leq C_{\alpha} 
		\bigg(\frac{1}{n}\bigg)^{\frac{1}{2}},
\end{equation}
\noindent  where $C_{\alpha}$ is a constant depending only on $\alpha$ and the uniform norms of $b$ and $\sigma$. Let $Y_n^g(t)=F_n(g)(t)-F(g)(t)$ and $Z_n^g(t)=|Y_n^g(t)|^2$. For any $0<\delta <e^{-1}$, define $\tau_n(g)=\inf\{ t\geq 0,|Y_n^g(t)|>\delta \}$. Given $\rho >0$, define
$$\Phi_{\rho}(\xi)=e^{\psi_{\rho}(\xi )},$$
where
$$\psi_{\rho}(\xi)=\int_0^{\xi}\frac{ds}{s\log{\frac{1}{s}}+\rho}.$$
Then
\begin{equation}
\Phi_{\rho}^{\prime}(\xi)(\xi\log{\frac{1}{\xi}}+\rho)=\Phi_{\rho}(\xi)
\end{equation}
By the chain rule,
$$\Phi_{\rho}(Z_n^g(t\wedge \tau_n(g)))$$
$$=1+2\int_0^{t\wedge \tau_n(g)}\Phi_{\rho}^{\prime}(Z_n^g(s))\bigg <Y_n^g(s),  b(F_n(g)(\frac{[ns]}{n})\big)-b(F(g)(s))\bigg >ds$$
\begin{equation}
+2\int_0^{t\wedge \tau_n(g)}\Phi_{\rho}^{\prime}(Z_n^g(s))\bigg <Y_n^g(s), ( \sigma (F_n(g)(\frac{[ns]}{n})\big)-\sigma (F(g)(s)))\dot{g}(s)\bigg >ds
\end{equation}
Using $(H.1)$ and (76), for $s\leq \tau_n(g)$,
$$|\bigg <Y_n^g(s),  b(F_n(g)\bigg(\frac{[ns]}{n}\bigg)\big)-b(F(g)(s))\bigg >|$$
$$\leq |\bigg <Y_n^g(s),  b(F_n(g)\bigg(\frac{[ns]}{n}\bigg)\big)-b(F_n(g)(s))\bigg >|$$
$$+|\bigg <Y_n^g(s),  b(F_n(g)(s))-b(F_n(g)(s))\bigg >|$$
$$\leq C |Y_n^g(s)| |F_n(g)\bigg(\frac{[ns]}{n}\bigg)-(F_n(g)(s)|\log (\frac{1}{|F_n(g)(\frac{[ns]}{n})-(F_n(g)(s)|})$$
$$+ \frac{1}{2} C Z_n^g(s) \log(\frac{1}{Z_n^g(s)})$$
$$\leq C C_{\alpha} (\frac{1}{n})^{\frac{1}{2}}\log (n^{\frac{1}{2}})+CZ_n^g(s) \log(\frac{1}{Z_n^g(s)})$$
\begin{equation}
\leq C(Z_n^g(s) \log(\frac{1}{Z_n^g(s)})+\rho)
\end{equation}
for $n\geq  N_{\alpha}^1$, where $N_{\alpha}^1$ depends  on $\alpha$ and $\rho$.
\noindent Similarly, for $s\leq \tau_n(g)$ and $n\geq N_{\alpha}^1$,
$$|Y_n^g(s)| ||\sigma (F_n(g)\bigg(\frac{[ns]}{n}\bigg)\big)-\sigma (F(g)(s))||$$
\begin{equation}
\leq C(Z_n^g(s) \log(\frac{1}{Z_n^g(s)})+\rho)
\end{equation}
It follows from (78) that for $n\geq  N_{\alpha}^1$ and all $g\in \{g; e(g)\leq \alpha \}$,
$$
\Phi_{\rho}(Z_n^g(t\wedge \tau_n(g)))$$
$$\leq 1+C \int_0^{t\wedge \tau_n(g)}\Phi_{\rho}^{\prime}(Z_n^g(s))(\frac{1}{Z_n^g(s)})+\rho)(1+|\dot{g}(s)|)ds$$
\begin{equation}
\leq C\int_0^t \Phi_{\rho}(Z_n^g(s\wedge \tau_n(g)))(1+|\dot{g}(s)|)ds
\end{equation}
By Gronwall's lemma,
$$\Phi_{\rho}(Z_n^g(1\wedge \tau_n(g)))\leq Ce^{1+\int_0^1|\dot{g}(s)|ds}$$
Since $\Phi_{\rho}(\xi)$ is incresing in $\xi$, it follows that for $n\geq  N_{\alpha}^1$,
\begin{equation}
\Phi_{\rho}(\sup_{ g\in \{g; e(g)\leq \alpha \}}Z_n^g(1\wedge \tau_n(g)))\leq Ce^{1+\alpha}
\end{equation}
Consequently, for any $\rho>0$,
\begin{equation}
\limsup_{n\rightarrow \infty}\Phi_{\rho}(\sup_{ g\in \{g; e(g)\leq \alpha \}}Z_n^g(1\wedge \tau_n(g)))\leq Ce^{1+\alpha}
\end{equation}
To complete the proof  it suffices  to show that for any $\delta>0$ there exists an integer $N$ such that if $n\geq N$, then $\tau_n(g)>1$ for all $g\in \{ g; e(g)\leq \alpha \}$. This is now a consequence of (83). In fact, otherwise, there exists $\delta>0$,  a subsequence $\{ n_k, k\geq 1\}$ of positive integers and $g_{n_k}\in \{g; e(g)\leq \alpha \}$ such that $\tau_{n_k}(g_{n_k})>1$. This implies that 
$$ 
\Phi_{\rho}(\sup_{ g\in \{g; e(g)\leq \alpha \}}Z_{n_{k}}^g(1\wedge \tau_{n_k}(g)))\geq \Phi_{\rho}(Z_{n_k}^{g_{n_k}}(1\wedge \tau_{n_k}(g_{n_k})))\geq \Phi_{\rho}(\delta^2)
$$
Combing this with (83), we get 
\begin{equation}
\Phi_{\rho}(\delta^2)\leq Ce^{1+\alpha}
\end{equation}
for all $\rho$. This leads to a contradiction since the left side of (84) tends to infinity as $\rho$ goes to $0$. The proof is complete.

\smallskip
\vskip 0.3cm

\noindent {\bf Proof of Theorem 4.1 when $b,\sigma$ are bounded }
 \medskip
Notice that $X_n^{\varepsilon}(s)$  
$=F_n(\varepsilon^{\frac{1}{2}}W)(s)$, where $W$ is the Brownian motion. 
The theorem follows from Proposition 5.2 , Lemma 5.3  and Theorem 4.2.23 in [DZ].

\section{Large deviations:  general case}

In this section  we  will  remove the boundeness assumptions on $b$ and $\sigma$. We begin with
\vskip 0.3cm

\noindent {\bf Proposition 6.1} 	\label{prop:threetwo}
{\it Assume $(H.2)$.
Then 
\begin{equation}  \label{eq:forty}
	\lim_{R \to \infty}\limsup_{\varepsilon \rightarrow 0} 
	\varepsilon \log P\big(\sup_{0\leq t\leq 1}
	\vert X^{\varepsilon}(t)-x_0\vert > R \big)=-\infty
\end{equation}
where  $X^{\varepsilon}(\cdot)$ is  the solution to equation (51).}
\vskip 0.3cm

\noindent {\bf Proof}. Let $\delta_0$ be a fixed small positive constant, say $\delta_0<\frac{1}{2}$.  Let $f\in C^1(R_{+})$ be a strictly positive  $C^1$ function on $R_{+}$ that satisfies
\[ f(s)=\left \{ \begin{array}{ll} -s\log s &\mbox{if $0\leq s\leq 1-\delta_0$}\\ s\log s & \mbox{ if $s\geq 1+\delta_0$}\end{array} \right.\]
From now on, we will use $C$ to denote a generic constant which may change from line to line.
It is easy to see that there exists a positive contant $C$ such that
\begin{equation}
s|\log s|+1\leq C(f(s)+1),\quad f^{\prime}(s)\geq -C  \quad \mbox{for $s\geq 0$}
\end{equation}
Define
$$ \psi(\xi) = \int_0^\xi \frac{ds}{f(s)+1}\quad \xi\geq 0. $$
and put
$$ \Phi_{\lambda}(\xi)=e^{\lambda\psi(\xi) }, \quad \lambda \geq 0$$ 
It follows from (86) that
\begin{equation}
 \Phi_{\lambda}^{\prime}(\xi)=\lambda \Phi_{\lambda}(\xi)\frac{1}{f(\xi)+1}\leq C\lambda \Phi_{\lambda}(\xi)\frac{1}{(\xi|\log \xi |+1)},
\end{equation}
and
$$
\Phi_{\lambda}^{\prime\prime }(\xi)=\lambda^2 \Phi_{\lambda}(\xi)\frac{1}{(f(\xi)+1)^2}-\lambda \Phi_{\lambda}(\xi)\frac{f^{\prime}(\xi)}{(f(\xi)+1)^2}
$$
\begin{equation}
\leq C(\lambda^2 +\lambda)\Phi_{\lambda}(\xi)\frac{1}{(\xi|\log \xi |+1)^2}
\end{equation}
Let $ \eta^{\varepsilon}(t)=X^{\varepsilon}(t)-x_o$ and $ \xi^{\varepsilon}(t)=|\eta^{\varepsilon}(t)|^2$. Define
$$\tau_R =\inf{\big\{t>0,  |\eta^{\varepsilon}(t)|\geq R\,\big\}},\quad R>0.$$
Now by It\^o formula, we have
\begin{eqnarray}
\Phi_{\lambda}(\xi^{\varepsilon}(t\wedge\tau_R))
&=& 1 + 2\varepsilon^{\frac{1}{2}} \int_0^{t\wedge\tau_R} \Phi_{\lambda}^{\prime}(\xi^{\varepsilon}(s))
\big <\eta^{\varepsilon}(s), \sigma(X^{\varepsilon}(s))dW_s\big >\nonumber \\
&+& 2\int_0^{t\wedge\tau_R} \Phi_{\lambda}^{\prime}(\xi^{\varepsilon}(s))\,
\big <\eta^{\varepsilon}(s), b(X^{\varepsilon}(s))\big >\,ds\nonumber \\
&+&\varepsilon \int_0^{t\wedge\tau_R}\Phi_{\lambda}^{\prime}(\xi^{\varepsilon}(s))
||\sigma(X^{\varepsilon}(s))||^2\,ds\nonumber \\
&+& 2 \varepsilon  \int_0^{t\wedge\tau_R} \Phi_{\lambda}^{\prime\prime }(\xi^{\varepsilon}(s))\,
|\sigma^*(X^{\varepsilon}(s))\eta^{\varepsilon}(s)|^2\,ds.
\end{eqnarray}
By $(H2)$, there exists $C_1>0$ such that
$$ \left\{\matrix{||\sigma(x)||^2&\leq& C_1\, \bigl(|x-x_o|^2\,\log{|x-x_o|}+1\bigr),\cr
|b(x)|&\leq& C_1\,  \bigl(|x-x_o|\,\log{|x-x_o|}+1\bigr).}\right.$$
It follows that
$$ {|\sigma^*(X^{\varepsilon}(s))\eta^{\varepsilon}(s)|^2\over (\xi^{\varepsilon}(s)\, |\log{\xi^{\varepsilon}(s)}|+1)^2 }
\leq C_1\ {\xi^{\varepsilon}(s)\, (\xi^{\varepsilon}(s)\, |\log{\xi^{\varepsilon}(s)}|+1)\over (\xi^{\varepsilon}(s)\, |\log{\xi^{\varepsilon}(s)}|+1)^2 }$$
which is dominated by a constant $C$. According to $(88)$, we get
\begin{equation}
 \int_0^{t\wedge\tau_R}\Phi_{\lambda}^{\prime\prime }(\xi^{\varepsilon}(s))\,|\sigma^*(X^{\varepsilon}(s))\eta^{\varepsilon}(s)|^2\,ds
\leq C(\lambda^2+\lambda) \int_0^{t\wedge\tau_R} \Phi_{\lambda}(\xi^{\varepsilon}(s))\, ds.
\end{equation}
In the same way, for some constant $C>0$, we have
\begin{equation}
 {|\big <\eta^{\varepsilon}(s), b(X^{\varepsilon}(s))\big >| + ||\sigma(X^{\varepsilon}(s))||^2\over
\xi^{\varepsilon}(s)\, |\log{\xi^{\varepsilon}(s)}| + 1}\leq C, \quad s>0.
\end{equation}
Combining above inequalities together, we get
$$E\big (\Phi_{\lambda}(\xi^{\varepsilon}(t\wedge\tau_R))\big )
\leq 1+ C (\varepsilon \lambda^2+\lambda)\int_0^t E\big (\Phi_{\lambda}(\xi^{\varepsilon}(t\wedge\tau_R))\big )\, ds,$$
which implies that
$$E\big (\Phi_{\lambda}(\xi^{\varepsilon}(1\wedge\tau_R))\big )
\leq e^{C (\varepsilon \lambda^2+\lambda)}.$$
Let $\lambda =\frac{1}{\varepsilon}$. It follows  that
$$ P\big(\sup_{0\leq t\leq 1}
	\vert X^{\varepsilon}(t)-x_0\vert > R \big) \Phi_{\frac{1}{\varepsilon}}(R)\leq E\big (\Phi_{\lambda}(\xi^{\varepsilon}(1\wedge\tau_R))\big )
$$
$$\leq e^{2C\frac{1}{\varepsilon}}$$
This gives that 
\begin{equation}
\limsup_{\varepsilon \rightarrow 0} 
	\varepsilon \log P\big(\sup_{0\leq t\leq 1}
	\vert X^{\varepsilon}(t)-x_0\vert > R \big)\leq C-\psi (R)
\end{equation}
Note that $\lim_{R\rightarrow \infty}\psi (R)=+\infty$. Letting $R$ tend to $+\infty$ in (92) proves the proposition.
\vskip 0.3cm 
	
\noindent For $R>0$, define $m_R=\sup\{ |b(x)|,||\sigma (x)||; |x|\leq R\}$ and 
$b_i^R=(-m_R -1)\vee b_i \wedge (m_R+1)$, $\sigma_{i,j}^R=(-m_R -1)\vee \sigma_{i,j} \wedge (m_R+1)$, $1\leq i\leq d$, $0\leq j\leq m$. Put
$b_R=(b_1^R,b_2^R,...,b_d^R)$ and $\sigma_R=(\sigma_{i,j}^R )_{1\leq i \leq d, 1\leq j\leq m}$. Then for $|x|\leq R$,
$$
	b_R(x)=b(x),\quad \sigma_R (x)=\sigma (x).
$$
\noindent and $b_R , \sigma_R$ satisfy  (H.1) and (H.2) with  the same  constants.

\noindent Let $X_R^{\varepsilon}(\cdot)$ be the solution to 
\begin{eqnarray}  \label{eq:fiftyfour}
	X_R^{\varepsilon}(t)&=&x_0 + \int_0^t b_R(X_R^{\varepsilon}(s))ds\nonumber \\
	  &+& \varepsilon^{\frac{1}{2}}\int_0^t
		\sigma_R (X_R^{\varepsilon}(s))dW_s, \,\, t > 0.
\end{eqnarray}

\noindent For $g$ with $e(g)<\infty$, let $F_R(g)$ be the solution to
\begin{eqnarray}	\label{eq:sixtythree}
	F_R(g)(t)&=&x_0+\int_0^t b_R\big(F_R(g)(s)\big)ds\nonumber \\
	  && + \int_0^t \sigma_R\big (F_R(g)(s)\big)\dot{g}(s)ds
\end{eqnarray}

\noindent Define 
\begin{equation}	\label{eq:sixtyfour}
  I_R(f)=\inf\bigg\{\frac{1}{2}e(g)^2; F_R(g)=f\bigg\}, \quad \quad 
	f\in C_{x_0}\big([0,1]\rightarrow R^d\big)
\end{equation}
\noindent If $\sup_{0\leq t\leq 1}|F(g)(t)|\leq R$, then $F(g)=F_R(g)$. Therefore, 
\begin{equation}	\label{eq:sixtyfive}
  I(f)=I_R(f), \quad \mbox{ for $f$ with } \sup_{0\leq t\leq 1}|f(t)|\leq R
\end{equation}

\noindent {\bf Lemma 6.2} 	\label{lem:threethree}
	{\it $I(\cdot )$ is a good rate function on $C_{x_0}([0,1], R^d)$,i.e., for any $\alpha \geq 0$,
the level set $\{f; I(f)\leq \alpha \}$ is compact.}
\vskip 0.3cm

\noindent {\bf Proof}. Arguing as in Lemma 5.3, it is easy to see that for $\alpha >0$, $\sup_{\{ g; e(g)\leq \alpha \}}||F(g)||_{\infty}\leq R$ for some constant $R$. Thus $F_R(\cdot )=F(\cdot )$ on $\{ g; e(g)\leq \alpha \}$. On the other hand,
It is easy to see that  $F_R(\cdot )$ is continuous on the  level set $\{g; e(g)\leq \alpha \}$, so is $F(\cdot )$. This  is sufficient  to conclude that $I(\cdot )$ is a good rate functional since $e(\cdot )$ is.
\vskip 0.3cm

\noindent {\bf Proof of Theorem 4.1 in the Unbounded Case.}

\vskip 0.3cm
\noindent Let $ \mu_{\varepsilon}^R$ denote the law of $X_R^{\varepsilon}(\cdot)$ on $C_{x_0}([0,1], R^d)$.  According to previous section, $\{\mu_{\varepsilon}^R, \varepsilon >0\}$ satisfies a large deviation principle with good rate function $I_R(\cdot)$. Note that $ \mu_{\varepsilon}^R$ and $ \mu_{\varepsilon}$ coincide on the ball $\{f; ||f||_{\infty}\leq R \}$. For $R>0$ and a closed subset $C\subset C_{x_0}([0,1], R^d)$, set $C_R=C\cap \{f; ||f||_{\infty}\leq R \}$.  Then,
\begin{eqnarray}	\label{eq:sixtysix}
	\mu_{\varepsilon} (C) &\leq & \mu_{\varepsilon}(C_{R})
		+P(\sup_{0\leq t\leq 1}|X^{\varepsilon}(t)|>R)\nonumber 	\\
	&= &\mu_{\varepsilon}^R(C_{R})+ P(\sup_{0\leq t\leq 1}|X^{\varepsilon}(t)|>R)
\end{eqnarray}
\noindent By  the large deviation principle for $\{\mu_{\varepsilon}^R, \varepsilon >0\}$,
$$
	\limsup_{\varepsilon \rightarrow 0}\, \varepsilon \log\mu_{\varepsilon} (C) 
$$
$$
\leq \big(-\inf_{f\in C_{R}} I_R(f)\big)
		\vee \big(\limsup_{\varepsilon \rightarrow 0} \varepsilon 
		\log P(\sup_{0\leq t\leq 1}|X^{\varepsilon}(t)|>R) \big)
$$
\begin{equation}
=\big(-\inf_{f\in C_{R}} I(f)\big)
		\vee \big(\limsup_{\varepsilon \rightarrow 0} \varepsilon 
		\log P(\sup_{0\leq t\leq 1}|X^{\varepsilon}(t)|>R) \big)
\end{equation}
\noindent Applying Proposition 6.1  and Letting  $R\rightarrow \infty$, we obtain
\begin{equation}	\label{eq:sixtynine}
	\limsup_{\varepsilon \rightarrow 0} \varepsilon \log\mu_{\varepsilon} 
		(C)\leq \big (-\inf_{f\in C} I(f)\big)
\end{equation}	
which is the upper bound.

\noindent Let $G$ be an open subset of $ C_{x_0}([0,1]\rightarrow R^d)$. Fix any $\phi_0 \in G$. Choosing $\delta>0$ such that $B(\phi_0,\delta)=\{f; ||f-\phi_0||_{\infty}\leq \delta \} \subset G$. Let $R=||\phi_0||_{\infty}+\delta$. Since 
$$
B(\phi_0,\delta)\subset \{f; ||f||_{\infty}\leq R\}
$$
We have 
\begin{eqnarray}	\label{eq:seventy}
	-I(\phi_0)=-I_R(\phi_0 )& \leq& \limsup_{\varepsilon \rightarrow 0} \, \varepsilon 
		\log\mu_{\varepsilon}^{R} \bigg(B(\phi_0,\delta)\bigg) \nonumber  \\
&=&\limsup_{\varepsilon \rightarrow 0} \, \varepsilon 
		\log\mu_{\varepsilon} \bigg(B(\phi_0,\delta)\bigg)\nonumber \\
	&\leq& \big(\limsup_{\varepsilon \rightarrow 0} \varepsilon 
		\log\mu_{\varepsilon}(G) \big)
\end{eqnarray}

\noindent Since $\phi_0$ is arbitrary, it follows that
\begin{equation}	\label{eq:seventytwo}
  -\inf_{f\in G}I(f )\leq \limsup_{\varepsilon \rightarrow 0} \varepsilon \log\mu_{\varepsilon}(G)
\end{equation}
which is the lower bound.

\vskip 0.3cm
\noindent {\bf Acknowledgements}. The second named author would like to thank I.M.B, UFR Sciences et techniques, Universit\'e de Bourgogne  for the support and the hospitality during  his stay where this work was initiated. This work is also partially supported by the British EPSRC , grant no. GR/R91144.

\end{document}